# Dynamical State Feedback Control for Linear Input Delay Systems, Part I: Dissipative Stabilization via Semidefinite Programming

Qian Feng, *Member, IEEE*, Cong Zhang, *Student Member, IEEE*, Bo Wei, *Member, IEEE*

***Abstract*— We propose an SDP-based framework to address the stabilization of input delay systems while taking into account dissipative constraints. A key to our approach is the introduction of the concept of parameterized linear dynamical state feedbacks (LDSFs), which draws inspiration from recent advancements in the analyses of distributed delays. The parameterized LDSFs generalize conventional predictor controllers, where the interpretation of state prediction is concealed and their degree of parameterization can be increased by adjusting the integral kernels. A sufficient condition for the existence of dissipative LDSFs is formulated as matrix inequalities by constructing a complete type Krasovskiĭ functional. To solve the bilinear matrix inequality in the synthesis condition, we employ an off-line inner convex approximation algorithm that can be initialized using the gains of predictor controllers obtained via explicit construction. So the unknowns of our LTDS can be computed by solving convex semidefinite programs. Numerical examples and simulations were experimented to demonstrate the validity and effectiveness of our methodology.**

## I. Introduction

It is well known [1]–[3] that input delay system (IDS) $\dot{\boldsymbol{x}}(t) = A\boldsymbol{x}(t) + B\boldsymbol{u}(t-r)$ is stabilized by predictor controller (PC)

$$\boldsymbol{u}(t) = K\boldsymbol{x}(t+r) = K\left(\mathsf{e}^{Ar}\boldsymbol{x}(t) + \int_{-r}^{0}\mathsf{e}^{-A\tau}B\boldsymbol{u}(t+\tau)\mathsf{d}\tau\right) \quad (1)$$

for any $r > 0$, if $\boldsymbol{x}(t)$ is measurable and $A + BK$ is Hurwitz. By utilizing PCs, stabilizing IDSs with $\boldsymbol{u}(t-r)$ is converted into an LTI state feedback problem, effectively circumventing the difficulties associated with the infinite-dimensional nature of input delays. State-of-the-art methods for PC design have been developed along several lines, including the algebraic/factorization approach [4]–[6], backstepping design using the PDE-ODE representations [7]–[9], and systems with special structures [10], [11].

However, exact state prediction is unrealizable for an IDS

$$\dot{\boldsymbol{x}}(t) = A\boldsymbol{x}(t) + B\boldsymbol{u}(t-r) + D\boldsymbol{w}(t) \quad (2)$$

affected by unmeasurable disturbances $\boldsymbol{w}(t)$, which requires a controller $\boldsymbol{u}(t) = K\mathsf{e}^{Ar}\boldsymbol{x}(t) + K\int_{-r}^{0}\mathsf{e}^{-A\tau}B\boldsymbol{u}(t+\tau)\mathsf{d}\tau + K\mathsf{e}^{Ar}\int_{0}^{r}\mathsf{e}^{-A\tau}D\boldsymbol{w}(t+\tau)\mathsf{d}\tau$ that depends on the future values of $\boldsymbol{w}(t)$. Several modified prediction schemes [12]–[15] for disturbance rejection have been developed, which require assumptions on the smoothness of $\boldsymbol{w}(t)$ and that $B = D$.

Alternatively, the presence of unmeasurable $\boldsymbol{w}(t)$ in (2) can be addressed from a different perspective. For practical systems, we can consider disturbance attenuation as the primary goal, which can be enforced by some performance objectives such as $\mathcal{H}^{\infty}$ control [5] or dissipative constraints [16], as explored in the field of robust control [17]. M. Meinsma [5], L. Mirkin [18] and G. Tadmor [19] with their co-authors made significant contributions to the robust control of systems with I/O delays in the frequency domain subject to various performance objectives such as $\mathcal{H}^{\infty}$ [5], [20], $\mathcal{H}^{2}$ [6], [21] and $\mathcal{L}^{1}$ [22]. In parallel, necessary and sufficient conditions for the existence of $\mathcal{H}^{\infty}$ predictor controllers have been derived via the Riccati equations in the time domain, see the works of G. Tadmor and A. Kojima in [23]–[27].

One potential issue with predictor controllers is that the numerical approximation of the integral equation using discrete delays may not be robust [28], [29] if $\left\|\int_{-r}^{0}K\mathsf{e}^{A\tau}B\mathsf{d}\tau\right\|$ is too "large" [30], as the hidden neutral terms [30, Eq. (12.7)] can destroy the stability of the approximated closed-loop system (CLS). This problem can be solved by adding a low-pass filter [28] to the discretized controller or by using a linear dynamical state feedback (LDSF) [1], [31] controller

$$\dot{\boldsymbol{u}}(t) = (KB + X)\,\boldsymbol{u}(t) + (KA - XK)\left(\mathsf{e}^{Ar}\boldsymbol{x}(t) + \int_{-r}^{0}\mathsf{e}^{-A\tau}B\boldsymbol{u}(t+\tau)\mathsf{d}\tau\right) \quad (3)$$

with any Hurwitz matrix $X \in \mathbb{R}^{p\times p}$, which is always constructible given the original PC in (1). Unlike (1), a CLS with (3) for any $r > 0$ is a retarded functional differential equation (FDE), whose properties ensure that the stability of CLS is robust to both pointwise delay implementations and tolerable delay mismatches [31]. However, these modifications, alter the controller dynamics such that the implemented controller with safe features $\widetilde{K}_N(s)$ does not converge to the nominal design of predictor controller $K(s)$ as the quadrature step $N \to \infty$ (i.e., $\lim_{N\to\infty}\widetilde{K}_N(s) \neq K(s)$). Yet, all ideal $\gamma_{\text{opt}}$ for $\mathcal{H}^{\infty}$ control computed by the methods in [24]–[27] are based on predictor controllers $K(s)$ instead of the implemented version $\widetilde{K}_N(s)$. Consequently, the performance level $\gamma_{\text{opt}}$ derived for the ideal $K(s)$ is not guaranteed after controller implementation. Moreover, alternative reset mechanism [32] for implementing predictor controllers suffers similar issues as $K(s) \neq Q(s)$, where the implemented $Q(s)$ controller renders the closed-loop system hybrid.

To resolve the above issue, it is more reasonable to consider a controller with the safe implementation feature like in (3) from the very beginning of our design to ensure performance consistency. For this structure, the implemented controller $\widetilde{T}_N(s)$ converges to the nominal design $T(s)$ as $N \to \infty$ [31], ensuring that the optimized performance index is preserved in implementation in an asymptotic sense. The controller proposed in this note is partially inspired by the structure in (3), in which the original input signal $\boldsymbol{u}(t)$ becomes a new state of the resulting CLS.

Given the recent advancements of efficient SDP algorithms [33], [34], it is of research interest to investigate the possibility of developing SDP solutions to the input delay problem with dissipative constraints, while keeping the delay-compensation benefits of the predictor structure. There are two potential issues with the structure of (3) in SDP settings with dissipative constraints (DCs). First, the presence of the distributed delay (DD) in (3) must be addressed explicitly, given that the resulting CLS with (3) is a linear distributed delay system. Second, the structure in (3) is restrictive, as the controller gain $K$ does not directly parameterize the functions in $(KA - XK)\,\mathsf{e}^{-A\tau}B$. Consequently, we prefer to consider a

A preliminary version of this paper was presented at the 17th IFAC Workshop on Time Delay Systems TDS, 2022. This work was partially supported by the National Natural Science Foundation of China under Grant Nos. 62303180, 62373150 and 62273145, as well as by the Fundamental Research Funds for Central Universities under Grant Nos. 2023MS032 and 2023YQ002, China. The authors are affiliated with the School of Control and Computer Engineering, North China Electric Power University, Beijing, 100096, China. Emails: qianfeng@ncepu.edu.cn, qfen204@aucklanduni.ac.nz, congzhang@ncepu.edu.cn, bowei@ncepu.edu.cn. Corresponding Author: Qian Feng more general form of controller parametrization that generalizes (3) to ensure improved performance. With the mathematical toolbox for DDs developed in [35], we can directly handle the functions in $(KA - XK)\,\mathsf{e}^{-A\tau}B$ via the Kronecker decomposition (KD). Consequently, we can now develop an SDP-based framework, via the Krasovskiĭ functional (KF) approach, to tackle input delay stabilization problems involving DDs while incorporating DCs in the controller design. A distinct advantage of the KF approach [36] lies in its inherent extensibility. Once established, extending the synthesis conditions to include various control ingredients, such as uncertainties and multi-objective criteria, is much easier than frequency domain approaches in general [30].

Assuming that all states $\boldsymbol{x}(t)$ are measurable, we propose an optimization framework for the stabilization of (2) with quadratic dissipative constraints based on the KF approach. The first distinguishing feature of our framework is that all FDEs in this study are defined in the extended sense using the Carathéodory interpretation [37, section 2.6] with respect to the Lebesgue measure. This formulation is particularly useful in modeling practical systems that are subject to noise and glitches. To address the presence of $\boldsymbol{u}(t-r)$ within an optimization framework, we introduce the concept of parameterized LDSFs that include both pointwise and DDs in $\boldsymbol{x}(t)$ and $\boldsymbol{u}(t)$. Parameterizing the DD integral is performed via the KD approach developed in [35] for the analysis of DDs in control contexts. From the perspective of functional analysis [38], the right-hand side of our LDSF can be interpreted as a special case of the parameterization of linear operators in the underlying Hilbert space, using particular sets of basis functions. In fact, our LDSF generalizes the structure in (3) with a much greater degree of parameterization, where the kernels of the DD integral always include the functions in $(KA - XK)\,\mathsf{e}^{-A\tau}B$. Moreover, we can add an arbitrarily many linearly independent functions to the DD kernel matrix, where these functions satisfy $\frac{\mathsf{d}\boldsymbol{f}(\tau)}{\mathsf{d}\tau} = M\boldsymbol{f}(\tau)$ for some $M \in \mathbb{R}^{d\times d}$ with $d \in \mathbb{N}$. This guarantees the parameters in (3) can always be selected as a particular choice of gains for the proposed controller that fully compensates the delay in $\boldsymbol{u}(t-r)$, while making sure that our LDSF has greater potential to fulfil DCs better than (3).

Expressed as a nonconvex semidefinite programming (SDP) problem, a sufficient condition for the existence of stabilizing LDSFs while satisfying the dissipativity condition is obtained by constructing a general KF, where the integral kernels of the KF coincide with the DD kernels in our LDSF. To solve the resulting bilinear matrix inequality (BMI) in the synthesis condition, an iterative algorithm is proposed based on the inner convex approximation [39], where each iteration is formulated as a convex SDP problem. The algorithm requires a feasible solution to start with, which can always be supplied by the values in (3) hanks to its connection with our LDSF and (1). Consequently, the algorithm does not rely on auxiliary convexification theorems that typically introduce additional conservatism [35]. This represents a key advantage of the proposed framework on account of the use of our controller that is related to both (3) and (1).

The remainder of this note is organized as follows. We formulate the synthesis problem using our parameterized LDSF controller and derive the corresponding CLS in Section II, explaining the structure of the proposed controller and its implications in control contexts. Then, the main framework for dissipative stabilization is presented in Section III, including a theorem and an algorithm. Finally, we show the effectiveness of our approach through a challenging example in Section IV before drawing concluding remarks in the final section.

**Notation:** Define $\mathbb{S}^n = \{X \in \mathbb{R}^{n\times n} : X = X^\top\}$. $\mathcal{C}(\mathcal{X};\mathbb{R}^n)$ stands for the Banach space of continuous functions with uniform norm $\sup_{\tau\in\mathcal{X}} \|\boldsymbol{f}(\tau)\|$, where $\|\mathbf{x}\|^2 = \mathbf{x}^\top\mathbf{x}$. Function space $\mathcal{M}(\mathcal{X};\mathbb{R}^d)$ contains all measurable functions from $\mathcal{X}$ to $\mathbb{R}^d$. Let $\mathcal{L}^p(\mathcal{X};\mathbb{R}^n) := \{\boldsymbol{f}(\cdot) \in \mathcal{M}(\mathcal{X};\mathbb{R}^n) : \|\boldsymbol{f}(\cdot)\|_p < +\infty\}$ with Lebesgue measurable set $\mathcal{X} \subseteq \mathbb{R}^n$ and seminorm $\|\boldsymbol{f}(\cdot)\|_p := \left(\int_{\mathcal{X}} \|\boldsymbol{f}(x)\|^p\,\mathsf{d}x\right)^{1/p}$. Define notation $\mathsf{Sy}(X) = X + X^\top$ for any $X \in \mathbb{R}^{n\times n}$. A row concatenation of objects/items is denoted by $\mathbf{Row}_{i=1}^n X_i = \left[X_1 \cdots X_i \cdots X_n\right]$. Symbol $*$ is applied to denote $[*]YX = X^\top YX$ or $X^\top Y[*] = X^\top YX$ or $\left[\begin{smallmatrix} A & B \\ * & C \end{smallmatrix}\right] = \left[\begin{smallmatrix} A & B \\ B^\top & C \end{smallmatrix}\right]$. Notation $\mathsf{O}_{n,m}$ stands for an $n \times m$ zero matrix that can be abbreviated as $\mathsf{O}_n$ with $n = m$, whereas $\mathbf{0}_n$ represents an $n\times 1$ zero column vector. We use symbol $\oplus$ to denote $X \oplus Y = \left[\begin{smallmatrix} X & \mathsf{O}_{n,q} \\ \mathsf{O}_{p,m} & Y \end{smallmatrix}\right]$ for any matrix $X \in \mathbb{R}^{n\times m}$, $Y \in \mathbb{R}^{p\times q}$ with its $n$-ary iterated form $\mathsf{diag}_{i=1}^{\nu} X_i = X_1 \oplus X_2 \oplus \cdots \oplus X_\nu$. Moreover, notation $\widetilde{\forall} x \in \mathcal{X}$ indicates for almost all $x \in \mathcal{X}$ with respect to the Lebesgue measure [38]. $\otimes$ is the Kronecker product. Finally, the precedents of matrix operations are *matrix (scalars) multiplications* $> \oplus > \otimes > +$.

# II. Problem formulation

The paper concerns the dissipative stabilization of IDS

$$
\begin{aligned}
\dot{\boldsymbol{x}}(t) &= A\boldsymbol{x}(t) + B\boldsymbol{u}(t-r) + D_1\boldsymbol{w}(t), \quad \widetilde{\forall} t \geq t_0 \in \mathbb{R} \\
\boldsymbol{z}(t) &= C_1\boldsymbol{\chi}(t) + C_2\boldsymbol{\chi}(t-r) \\
&\quad + \int_{-r}^{0} \widetilde{C}_3(\tau)\boldsymbol{\chi}(t+\tau)\mathsf{d}\tau + D_3\boldsymbol{w}(t), \quad r > 0
\end{aligned}
\tag{4}
$$

where $\boldsymbol{x}(t) \in \mathbb{R}^n$ is the state variable, $\boldsymbol{u}(t) \in \mathbb{R}^p$ stands for system's input, and $\boldsymbol{z}(t) \in \mathbb{R}^m$ indicates the regulated output with disturbances $\boldsymbol{w}(\cdot) \in \mathcal{L}^2\left([t_0, +\infty]; \mathbb{R}^q\right)$. We use $\boldsymbol{\chi}^\top(t) := \left[\boldsymbol{x}^\top(t) \;\; \boldsymbol{u}^\top(t)\right]$ to denote the augmented vector composed of both the state and input. The dimensions of $A \in \mathbb{R}^{n\times n}$, $B \in \mathbb{R}^{n\times p}$, $D_1 \in \mathbb{R}^{n\times q}$ and $C_1, C_2, \widetilde{C}_3(\tau) \in \mathbb{R}^{m\times\nu}$, $D_3 \in \mathbb{R}^{m\times q}$ are determined by $m, n, p, q \in \mathbb{N}$ with $\nu := n + p$. Finally, we assume that there exists $K \in \mathbb{R}^{p\times n}$ such that $A + BK$ is Hurwitz.

It is crucial to stress that the FDE in (4) is defined in the extended sense for almost all $t \geq t_0$ w.r.t the Lebesgue measure, even when $\boldsymbol{w}(t) \equiv \mathbf{0}_q$. Since the right-hand side of (4) is linear, the existence and uniqueness of the solution to (4) are ensured by the Carathéodory's existence and uniqueness theorem for FDEs [37, section 2.6]. FDEs in the extended sense are particularly useful in modeling engineering and cybernetic systems that are subject to noise and glitches, the dynamics of which are difficult to capture with conventional FDEs [37, Chapter 1].

## A. Notes on the Properties of Controller (3)

Assuming that the values of $\boldsymbol{x}(t)$ are measurable, the PC in (1) can always exponentially stabilize the nominal system in (4) with $\boldsymbol{w}(t) \equiv \mathbf{0}_q$, if there exists $K \in \mathbb{R}^{p\times n}$ such that $A + BK$ is Hurwitz. To ensure safe numerical implementation [28] for integral terms, one can employ the dynamical controller in (3), which is always applicable to the nominal IDS in (4) using any Hurwitz matrix $X \in \mathbb{R}^{p\times p}$. Indeed, the spectrum [31] of the nominal system in (4) with (3) is

$$
\Big\{ s \in \mathbb{C} : \det\left(sI_n - A - BK\right)\det\left(sI_p - X\right) = 0 \Big\},
\tag{5}
$$

which implies that exponential stability is always achievable for the CLS with some $K \in \mathbb{R}^{p\times n}$ for any Hurwitz $X \in \mathbb{R}^{p\times p}$.

## B. Concept of KD for the Analysis of DDs

To handle DDs explicitly within an optimization framework, we employ the Kronecker decomposition (KD) for the matrix-valued functions in this note. The idea of KD was originally developed in [35] to address the presence of nontrivial DDs in FDEs.

*Assumption 1:* There exist matrices $C_3 \in \mathbb{R}^{m\times d\nu}$, $\Gamma \in \mathbb{R}^{p\times d\nu}$ and function $\boldsymbol{f}(\cdot) \in C^1\left([-r, 0]; \mathbb{R}^d\right)$ such that $\widetilde{C}_3(\tau) = C_3F(\tau)$ and $\left[\mathsf{O}_{p,n}\ \ (KA{-}XK)\mathsf{e}^{-A\tau}B\right] = \Gamma F(\tau) \in \mathbb{R}^{p\times\nu}$, where $F(\tau) = \left(\sqrt{\mathfrak{F}^{-1}}\boldsymbol{f}(\tau)\otimes I_\nu\right)$ and function $\boldsymbol{f}(\cdot)$ satisfies

$$\begin{aligned} &\exists M \in \mathbb{R}^{d\times d}, \frac{\mathsf{d}\boldsymbol{f}(\tau)}{\mathsf{d}\tau} = M\boldsymbol{f}(\tau), \\ &\mathfrak{F} = \int_{-r}^{0}\boldsymbol{f}(\tau)\boldsymbol{f}^{\top}(\tau)\mathsf{d}\tau \succ 0, \end{aligned} \tag{6}$$

and $X \in \mathbb{R}^{p\times p}$ is any Hurwitz matrix with any $K \in \mathbb{R}^{p\times n}$ such that $A + BK$ Hurwitz. Moreover, $\frac{\mathsf{d}}{\mathsf{d}\tau}$ in (6) represents weak derivatives.

A key ingredient of the KD is the use of the Gramian matrix $\mathfrak{F}$ (see [40, Th. 7.2.10]) of vector-valued function $\boldsymbol{f}(\cdot)$ in (6). Inequality $\mathfrak{F} \succ 0$ indicates that the functions in each row of $\boldsymbol{f}(\cdot)$ are linearly independent in a Lebesgue sense over $[-r, 0]$. Note that there are no contradictions between $\operatorname{rank}\mathfrak{F} = d$ and the fact that $\operatorname{rank}\left[\boldsymbol{f}(\tau)\boldsymbol{f}^{\top}(\tau)\right] = 1$, as $\operatorname{rank}\mathfrak{F} \neq \operatorname{rank}\left[\boldsymbol{f}(\tau)\boldsymbol{f}^{\top}(\tau)\right] = 1$ due to the integral operation. A simple example for $\mathfrak{F}$ is the Hilbert matrices [40] such as $\int_0^1[1\ \ \tau]^{\top}[1\ \ \tau]\mathsf{d}\tau = \int_0^1\begin{bmatrix}1 & \tau\\ * & \tau^2\end{bmatrix}\mathsf{d}\tau = \begin{bmatrix}1 & 1/2\\ * & 1/3\end{bmatrix} \succ 0$.

The rationale behind the conditions in Assumption 1 is fully explained **later in this note** by considering the structure of our proposed controller and the procedures to derive our synthesis condition. In terms of mathematics, the conditions in Assumption 1 indicate that all functions in $\widetilde{C}_3(\tau)$ and the DD kernels $(KA - XK)\,\mathsf{e}^{-A\tau}B$ in (3) are included in $\boldsymbol{f}(\cdot)$. This regulates what functions can be in $\widetilde{C}_3(\tau)$ but does not impose any structural constraint on $A$ or $B$, as $\frac{\mathsf{d}\mathsf{e}^{-A\tau}B}{\mathsf{d}\tau} = -A\mathsf{e}^{-A\tau}B$ can always be rewritten in the form of the differential equation in (6). The functions in $\boldsymbol{f}(\cdot)$ are solutions to linear homogeneous differential equations with constant coefficients such as polynomials, trigonometric functions, exponential functions and their products. Moreover, it is noteworthy that we can always add an unlimited number of functions to $\boldsymbol{f}(\cdot)$ satisfying (6), since the dimensions of $\boldsymbol{f}(\tau)$ and $M$ are unbounded and new zeros can be added to $C_3$ and $\Gamma$ to ensure consistent decompositions. Finally, we emphasize that the use of $\sqrt{\mathfrak{F}^{-1}}$ does not affect the existence of $C_3$ and $\Gamma$ in Assumption 1 as $\sqrt{\mathfrak{F}^{-1}}$ is of full rank.

### *C. Concept of Parameterized LDSFs in Relation to (3)*

Inspired by the structure in (3) and our previous works on DDs in [35], we introduce the concept of parameterized LDSFs

$$\begin{aligned} \dot{\boldsymbol{u}}(t) = \mathfrak{c}\left(\boldsymbol{\chi}_t(\cdot)\right) &= K_1\boldsymbol{\chi}(t) + K_2\boldsymbol{\chi}(t-r) \\ &+ \textstyle\int_{-r}^{0}K_3F(\tau)\boldsymbol{\chi}(t+\tau)\mathsf{d}\tau, \quad \boldsymbol{\chi}_t(\theta) := \boldsymbol{\chi}(t+\theta), \end{aligned} \tag{7}$$

with unknowns $K_1; K_2 \in \mathbb{R}^{p\times\nu}$ and $K_3 \in \mathbb{R}^{p\times d\nu}$, where $F(\tau)$, $\boldsymbol{f}(\tau)$ and $\mathfrak{F}$ are given by Assumption 1. The structure in (7) is a parameterization of linear operator $\mathcal{C}(\mathcal{X};\mathbb{R}^n) \ni \boldsymbol{\psi}(\cdot) \mapsto \mathfrak{c}\left(\boldsymbol{\psi}(\cdot)\right) \in \mathbb{R}^p$ using $\boldsymbol{f}(\cdot)$. Note that all values of $\boldsymbol{\chi}(t)$ can be measured, as $\boldsymbol{x}(t)$ is measurable and the values of $\boldsymbol{u}(t)$ are always available.

Since all functions in $(KA - XK)\,\mathsf{e}^{-A\tau}B$ of (3) are covered by some functions in $\boldsymbol{f}(\cdot)$ by Assumption 1, the expression in (7) generalizes that of (3), as we can always select

$$\begin{aligned} &K_1 = \left[(KA - XK)\,\mathsf{e}^{Ar}\quad KB + X\right] \in \mathbb{R}^{p\times\nu}, \\ &K_2 = \mathsf{O}_{p,\nu}, \quad K_3 = \Gamma \in \mathbb{R}^{p\times d\nu}. \end{aligned} \tag{8}$$

This clearly explains the rationale of $\left[\mathsf{O}_{p,n}\ \ (KA{-}XK)\mathsf{e}^{-A\tau}B\right] = \Gamma F(\tau)$ in Assumption 1, which ensures that (3) is a special case of (7). As all terms in (7) are directly parameterized, controller (7) has greater potential to fulfil DCs than (3).

Using LDSFs to stabilize IDSs is not a new idea, as this concept was explored in [1], [28], [31] in various control contexts. However, existing studies on LDSFs assume that they are directly constructed from PCs, thereby inheriting the parameter structure of the underlying PCs as in (3). In contrast, the unknowns in (7) are not restricted by the $K$ in (1). In fact, interpreting the structure of (7) for state prediction (finite spectrum assignment), as is the case for a PC, is not necessary. Consequently, the concept specified in (7) offers a novel perspective on addressing the stabilization of IDSs when DCs are considered.

Since feedback loops can also be affected by external disturbances introduced by the operating environment, we use equation

$$\dot{\boldsymbol{u}}(t) = \mathfrak{c}\left(\boldsymbol{\chi}_t(\cdot)\right) + D_2\boldsymbol{w}(t), \tag{9}$$

to denote such cases as illustrated in Figure 1, where $D_2 \in \mathbb{R}^{p\times q}$ is known. Note that unknown disturbance $D_2\boldsymbol{w}(t)$ is **not** a part of controller (7). Once the unknowns in (7) are obtained, we employ controller (7) to open-loop system (4) instead of implementing (9). Utilizing (9) in the subsequent analysis, however, ensured that controller (7) is robust to disturbances with respect to the DCs.

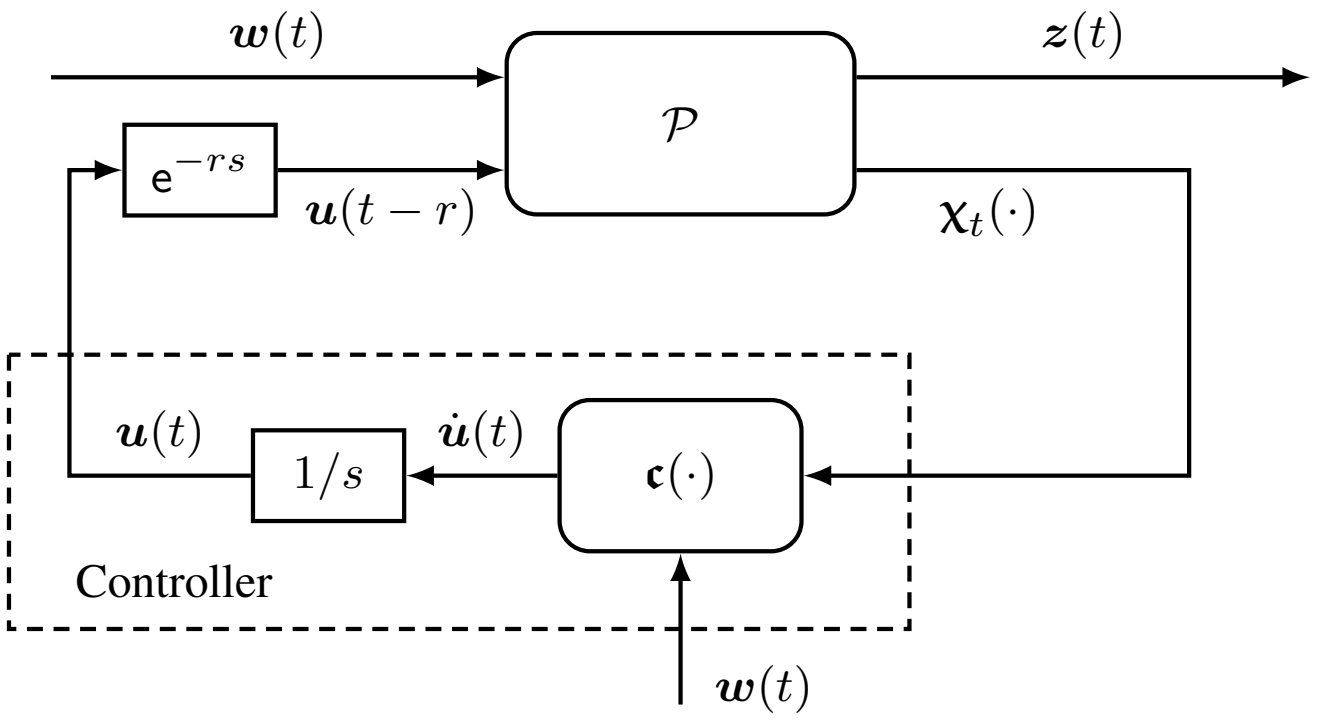

Fig. 1**:** Diagram and Configuration of the Closed-Loop Systems

### *D. Expression of Closed-Loop System*

By Assumption 1 and open-loop system (4) and controller (9), the dynamics of the CLS denoted in Figure 1 is

$$\begin{aligned} \dot{\boldsymbol{\chi}}(t) &= \begin{bmatrix}\begin{matrix}A & \mathsf{O}_{n,p}\end{matrix}\\ K_1\end{bmatrix}\boldsymbol{\chi}(t) + \begin{bmatrix}\begin{matrix}\mathsf{O}_n & B\end{matrix}\\ K_2\end{bmatrix}\boldsymbol{\chi}(t-r), \\ &+ \int_{-r}^{0}\begin{bmatrix}\mathsf{O}_{n,d\nu}\\ K_3\end{bmatrix}F(\tau)\boldsymbol{\chi}(t+\tau)\mathsf{d}\tau + \begin{bmatrix}D_1\\ D_2\end{bmatrix}\boldsymbol{w}(t), \\ \boldsymbol{z}(t) &= \boldsymbol{\Sigma}\boldsymbol{\vartheta}(t), \quad \forall\theta \in [-r, 0],\ \boldsymbol{\chi}(t_0+\theta) = \boldsymbol{\phi}(\theta), \\ \boldsymbol{\Sigma} &= \left[C_1\quad C_2\quad C_3\quad D_3\right],\ \boldsymbol{\chi}(t) := \left[\boldsymbol{x}^{\top}(t)\ \ \boldsymbol{u}^{\top}(t)\right]^{\top} \\ \boldsymbol{\vartheta}^{\top}(t) &= \left[\boldsymbol{\chi}^{\top}(t)\quad \boldsymbol{\chi}^{\top}(t-r)\quad \textstyle\int_{-r}^{0}\boldsymbol{\chi}^{\top}(t+\tau)F^{\top}(\tau)\mathsf{d}\tau\quad \boldsymbol{w}^{\top}(t)\right] \end{aligned} \tag{10}$$

$\widetilde{\forall} t \geq t_0 \in \mathbb{R}$, where $\boldsymbol{\phi}(\cdot) \in \mathcal{C}([-r, 0]; \mathbb{R}^{\nu})$ is the initial condition. Similar to (4), the FDE in (10) is defined in an extended sense w.r.t the Lebesgue measure using the Carathéodory framework. Given the structures in (10) and (7), the following conclusions are established.

*Conclusion 1:* The equation in (10) is a retarded FDE, which satisfies [31, Th. 2]. Namely, if (7) exponentially stabilizes (4) with $\boldsymbol{w}(t) \equiv 0$, then the DD in (7) can always be numerically implemented by ordinary quadrature rules without causing instabilities to the CLS, provided that the numerical accuracy is sufficient. Also, the stability of CLS is robust w.r.t to small delay mismatches [30, Eq. (12.7)].

*Conclusion 2:* Since (3) is a special case of (7), we can always utilize (8) for (7) to exponentially stabilize (4) with $\boldsymbol{w}(t) \equiv 0$ for any $r > 0$. Consequently, our LDSF (7) always has the ability to fully eliminate the effects of input delays, similar to (3) and (1).

## III. Dissipative Dynamical State Feedback

We first present the lemma to be used in determining the stability of the trivial solution of (10). It is crucial to point out that Lemma

1 is **not** a special case of the classical Krasovskiĭ stability theorem in [37, Section 5.2] for FDEs with classical derivatives, as the FDE in (10) is defined in the extended sense w.r.t Lebesgue measure.

*Lemma 1:* The trivial solution of the FDE in (10) with $\boldsymbol{w}(t) \equiv \mathbf{0}_q$ is globally exponentially stable, if there exist $\epsilon_1; \epsilon_2; \epsilon_3 > 0$ and a continuous functional $\mathsf{v} : \mathscr{C}([-r, 0]; \mathbb{R}^\nu) \to \mathbb{R}$ such that

$$\epsilon_1 \|\boldsymbol{\phi}(0)\|^2 \leq \mathsf{v}(\boldsymbol{\phi}(\cdot)) \leq \epsilon_2 \|\boldsymbol{\phi}(\cdot)\|_\infty^2 \tag{11a}$$

$$\widetilde{\forall} t \geq t_0, \ \tfrac{\mathsf{d}}{\mathsf{d}t}\mathsf{v}(\boldsymbol{\chi}_t(\cdot)) \leq -\epsilon_3 \|\boldsymbol{\chi}(t)\|^2 \tag{11b}$$

for any function $\boldsymbol{\phi}(\cdot) \in \mathscr{C}([-r, 0]; \mathbb{R}^\nu)$ in (10), where $\|\boldsymbol{\phi}(\cdot)\|_\infty^2 := \sup_{-r \leq \tau \leq 0} \|\boldsymbol{\phi}(\tau)\|^2$. Furthermore, operator $\boldsymbol{\chi}_t(\cdot)$ in (11b) is defined by $\forall t \geq t_0$, $\forall \theta \in [-r, 0]$, $\boldsymbol{\chi}_t(\theta) = \boldsymbol{\chi}(t + \theta)$ in which function $[t_0 - r, \infty) \ni t \mapsto \boldsymbol{\chi}(t) \in \mathbb{R}^\nu$ satisfies (10) with $\boldsymbol{w}(t) \equiv \mathbf{0}_q$.

*Proof:* The proof is similar to that of [41, Cor. 1], as the right-hand side of the FDE in (10) is linear and retarded, which also implies that exponential stability is inferred from asymptotic stability on the basis of the system's spectrum. ■

The following characterization of dissipativity is grounded in the original definition in [42].

*Definition 1:* CLS (10) with supply rate function (SRF) $\mathsf{s}(\boldsymbol{z}(t), \boldsymbol{w}(t))$ is said to be dissipative if there exists a continuous functional $\mathsf{v} : \mathscr{C}([-r, 0]; \mathbb{R}^\nu) \to \mathbb{R}$ such that

$$\widetilde{\forall} t \geq t_0, \quad \tfrac{\mathsf{d}}{\mathsf{d}t}\mathsf{v}(\boldsymbol{\chi}_t(\cdot)) - \mathsf{s}(\boldsymbol{z}(t), \boldsymbol{w}(t)) \leq 0 \tag{12}$$

with $t_0 \in \mathbb{R}$, $\boldsymbol{z}(t)$ and $\boldsymbol{w}(t)$ in (10). Moreover, the operator $\boldsymbol{\chi}_t(\cdot)$ in (12) is defined by $\forall t \geq t_0$, $\forall \theta \in [-r, 0]$, $\boldsymbol{\chi}_t(\theta) = \boldsymbol{\chi}(t + \theta)$ with function $[t_0 - r, \infty) \ni t \mapsto \boldsymbol{\chi}(t) \in \mathbb{R}^\nu$ satisfying (10).

To characterize dissipativity in this note, we select

$$\begin{gathered}\mathsf{s}(\boldsymbol{z}(t), \boldsymbol{w}(t)) = \begin{bmatrix} \boldsymbol{z}(t) \\ \boldsymbol{w}(t) \end{bmatrix}^\top \begin{bmatrix} \widetilde{J}^\top J_1^{-1} \widetilde{J} & J_2 \\ * & J_3 \end{bmatrix} \begin{bmatrix} \boldsymbol{z}(t) \\ \boldsymbol{w}(t) \end{bmatrix}, \\ \widetilde{J}^\top J_1^{-1} \widetilde{J} \preceq 0, \qquad \widetilde{J} \in \mathbb{R}^{m \times m}, \qquad \mathbb{S}^m \ni J_1^{-1} \prec 0, \\ J_2 \in \mathbb{R}^{m \times q}, \quad J_3 \in \mathbb{S}^q \end{gathered} \tag{13}$$

as the SRF, whose structure is grounded in the quadratic constraints in [17] with a slightly modified formulation. It can feature numerous performance objectives [17] such as

- $\mathscr{L}^2/\mathscr{H}^\infty$ Gain Performance: $J_1 = -\gamma I_m$, $\widetilde{J} = I_m$, $J_2 = \mathsf{O}_{m,q}$, $J_3 = \gamma I_q$ with $\gamma > 0$
- Strict Passivity: $J_1 \prec 0$, $\widetilde{J} = \mathsf{O}_m$, $J_2 = I_m$, $J_3 = \mathsf{O}_m$, $m = q$
- Numerous Sector Constraints in [16, Table 1].

With all the tools above, the main results of this note are presented in the following theorem, whose proof is provided in Appendix B.

*Theorem 1:* Let all the parameters in Assumption 1 be given. Then the CLS in (10) with the SRF in (13) is dissipative, and the trivial solution of (10) with $\boldsymbol{w}(t) \equiv \mathbf{0}_q$ is exponentially stable if there exist $P \in \mathbb{S}^\nu$, $Q \in \mathbb{R}^{\nu \times d\nu}$, $R \in \mathbb{S}^{d\nu}$, $S; U \in \mathbb{S}^\nu$ and $K_1 \in \mathbb{R}^{p \times \nu}$, $K_2 \in \mathbb{R}^{p \times \nu}$ and $K_3 \in \mathbb{R}^{p \times d\nu}$ such that

$$\begin{bmatrix} P & Q \\ * & R + I_d \otimes S \end{bmatrix} \succ 0, \quad S \succ 0, \quad U \succ 0, \tag{14a}$$

$$\boldsymbol{\Phi} + \mathsf{Sy}\left(\mathbf{P}^\top \begin{bmatrix} \boldsymbol{\Pi} & \mathsf{O}_{\nu,m} \end{bmatrix}\right) \prec 0, \tag{14b}$$

where $\boldsymbol{\Pi} = \mathbf{A} + \mathbf{B}\mathbf{K}$ and $\mathbf{P} = \begin{bmatrix} P & \mathsf{O}_\nu & Q & \mathsf{O}_{\nu,(q+m)} \end{bmatrix}$ and

$$\mathbf{A} = \begin{bmatrix} A & \mathsf{O}_{n,\nu} & B & \mathsf{O}_{n,d\nu} & D_1 \\ \mathsf{O}_{p,n} & \mathsf{O}_{p,\nu} & \mathsf{O}_p & \mathsf{O}_{p,d\nu} & D_2 \end{bmatrix}, \tag{15a}$$

$$\mathbf{B} = \begin{bmatrix} \mathsf{O}_{n,p} \\ I_p \end{bmatrix}, \ \mathbf{K} = \begin{bmatrix} K_1 & K_2 & K_3 & \mathsf{O}_{p,q} \end{bmatrix}, \tag{15b}$$

$$\begin{aligned}\boldsymbol{\Phi} = {} & \mathsf{Sy}\left(\begin{bmatrix} Q \\ \mathsf{O}_{\nu,d\nu} \\ R \\ \mathsf{O}_{(q+m),d\nu} \end{bmatrix} \begin{bmatrix} \mathbf{F} & \mathsf{O}_{d\nu,m} \end{bmatrix}\right) \\ & + (S + rU) \oplus [-S] \oplus [-I_d \otimes U] \oplus J_3 \oplus J_1^{-1} \\ & + \mathsf{Sy}\left(\begin{bmatrix} \mathsf{O}_{m,(2\nu+d\nu)} & J_2 & I_m \end{bmatrix}^\top \begin{bmatrix} \boldsymbol{\Sigma} & \mathsf{O}_m \end{bmatrix}\right)\end{aligned} \tag{15c}$$

$$\mathbf{F} = \begin{bmatrix} F(0) & -F(-r) & -\sqrt{\mathfrak{F}^{-1}} M \sqrt{\mathfrak{F}} \otimes I_\nu & \mathsf{O}_{d\nu,q} \end{bmatrix} \tag{15d}$$

with $\boldsymbol{\Sigma}$ in (10) and $F(\tau)$, $\mathfrak{F}$ in Assumption 1. The total number of unknowns is $(0.5d^2 + d + 1.5)\nu^2 + (0.5d + pd + 2p + 1.5)\nu \in \mathcal{O}(d^2\nu^2)$.

*Remark 1:* By using the techniques [36] in matrix analysis in light of the SDP constraint developed in (14b), the synthesis condition in Theorem 1 can be restructured to optimize the generalized $\mathscr{H}^2$ objective with $D_3 = \mathsf{O}_{m,q}$ based on its definition in [17], [43].

*Remark 2:* The rationale behind the condition in (6) is to ensure that both $\int_{-r}^{0} F(\tau)\boldsymbol{\chi}(t + \tau)\mathsf{d}\tau$ in (29) and its derivative can be expressed as products between constant matrices and $\boldsymbol{\vartheta}(t)$ in (10), as shown in (31) and (32). Moreover, $\mathfrak{F}$ and $F(\tau) = \sqrt{\mathfrak{F}^{-1}}\boldsymbol{f}(\tau) \otimes I_\nu$ in (6) and (29) are defined such that inequality (28) is always applicable in steps (33),(36),(37) with $Y = I_d$. This ensures that the numerical solvability of the SDP constraints in (14a)–(14b) is not affected by the condition number of $\mathfrak{F}$, since no $\mathfrak{F}$-related terms are found in the diagonal matrix blocks in the inequalities in (14a)–(14b).

The conservatism of Theorem 1 can be qualitatively estimated in relation to existing results. When $K_1; K_2; K_3$ are known, Theorem 1 becomes a stability analysis condition for (10), which amounts to the stability condition developed in [35]. The results in [35], in turn, generalize the stability conditions in [44], [45] wherein $\boldsymbol{f}(\tau)$ comprises a list of Legendre polynomials $\boldsymbol{\ell}_d(\tau)$ up to degree $d$. Notably, recent findings in [45] reveal that the resulting LMI conditions are both sufficient and necessary, given a sufficiently large $d$ for $\boldsymbol{\ell}_d(\tau)$, when DDs are absent in the system. In fact, our KF (29) can be interpreted as an approximation (parameterization) of the complete KF for $\dot{\boldsymbol{\chi}}(t) = A_0\boldsymbol{\chi}(t) + A_1\boldsymbol{\chi}(t - r)$ using basis function $\boldsymbol{f}(\cdot)$ rather than relying solely on $\boldsymbol{\ell}_d(\tau)$, where $\dim(\boldsymbol{f}(\tau))$ is also unbounded since an unlimited number of new functions can be added to $\boldsymbol{f}(\cdot)$ satisfying (6). Given also the robust testing results of the approaches in [35], [44], [45] for various challenging numerical examples, these facts suggest that Theorem 1 is nonconservative.

Inequality (14b) is bilinear as it contains $\mathbf{P}^\top\mathbf{B}\mathbf{K}$ and its transpose. The convexification strategies in [35], [36] can be applied to (14b) at the expense of introducing extra variables and potential conservatism. However, the problem we are addressing has a unique feature due to the expression in (7). Since we can always find some $K, X, \Gamma$ in Assumption 1 via (3), the parameters in (8) can be selected as a candidate for the unknowns in (9) that can render (14b) feasible. If the parameters in (8) are utilized and are known, constraint (14b) becomes convex, which can be computed by standard SDP solvers [33]. This property is particularly useful when applying iterative algorithms to (14b) aimed at finding locally optimal solutions, which always requires a feasible solution to (14b) as a starting point.

### *A. An Iterative Algorithm (offline) for the BMI in Theorem 1*

Given our useful findings in (8), we employ the concept in [39], which generalizes the convex-concave decomposition algorithm in [46], to construct an iterative algorithm (offline) for computing locally optimal solutions for the BMI in (14b). Consider function

$$\begin{aligned}\mathbb{S}^\ell \ni \Delta\left(\mathbf{G}, \widetilde{\mathbf{G}}, \mathbf{N}, \widetilde{\mathbf{N}}\right) := [*] \left[Z \oplus (I_\nu - Z)\right]^{-1} \begin{bmatrix} \mathbf{G} - \widetilde{\mathbf{G}} \\ \mathbf{N} - \widetilde{\mathbf{N}} \end{bmatrix} \\ + \mathsf{Sy}\left(\widetilde{\mathbf{G}}^\top\mathbf{N} + \mathbf{G}^\top\widetilde{\mathbf{N}} - \widetilde{\mathbf{G}}^\top\widetilde{\mathbf{N}}\right) + \mathbf{T}\end{aligned} \tag{16}$$

with $Z \oplus (I_\nu - Z) \succ 0$ satisfying $\forall \mathbf{G}; \widetilde{\mathbf{G}} \in \mathbb{R}^{\nu \times \ell}$, $\forall \mathbf{N}; \widetilde{\mathbf{N}} \in \mathbb{R}^{\nu \times \ell}$ :

$$\mathbf{T} + \mathsf{Sy}\left(\mathbf{G}^\top\mathbf{N}\right) = \Delta\left(\mathbf{G}, \mathbf{G}, \mathbf{N}, \mathbf{N}\right) \preceq \Delta\left(\mathbf{G}, \widetilde{\mathbf{G}}, \mathbf{N}, \widetilde{\mathbf{N}}\right). \tag{17}$$

This shows that $\Delta(\mathbf{G},\widetilde{\mathbf{G}},\mathbf{N},\widetilde{\mathbf{N}})$ is a PSD-overestimate [39, Definition 2.2] of $\mathbf{T}+\mathsf{Sy}\left[\mathbf{G}^\top\mathbf{N}\right]$ with respect to the parameterization

$$\mathbf{vec}(\widetilde{\mathbf{G}})=\mathbf{vec}(\mathbf{G}) \quad \& \quad \mathbf{vec}(\widetilde{\mathbf{N}})=\mathbf{vec}(\mathbf{N}),$$

where $\mathbf{vec}(\mathbf{G})$ is the vectorization [40] of $\mathbf{G}$. We rewrite (14b) as

$$\boldsymbol{\Theta}=\boldsymbol{\Phi}+\mathsf{Sy}\left(\mathbf{P}^\top\begin{bmatrix}\boldsymbol{\Pi} & \mathsf{O}_{\nu,m}\end{bmatrix}\right)=\widehat{\boldsymbol{\Phi}}+\mathsf{Sy}\left(\mathbf{P}^\top\mathbf{B}\boldsymbol{\mathfrak{K}}\right)\prec 0, \quad (18)$$

where $\boldsymbol{\mathfrak{K}}:=\begin{bmatrix}\mathbf{K} & \mathsf{O}_{p,m}\end{bmatrix}$ and $\widehat{\boldsymbol{\Phi}}:=\mathsf{Sy}\left(\mathbf{P}^\top\begin{bmatrix}\mathbf{A} & \mathsf{O}_{\nu,m}\end{bmatrix}\right)+\boldsymbol{\Phi}$ with $\mathbf{A},\mathbf{K}$ in (15a)–(15b). Note that the unknowns in $\widehat{\boldsymbol{\Phi}}$ are convex.

Now let $\mathbf{T}=\widehat{\boldsymbol{\Phi}}$, $\mathbf{G}=\mathbf{P}$ with matrices in Theorem 1, and

$$\begin{aligned}
&\widetilde{\mathbf{G}}=\widetilde{\mathbf{P}}:=\begin{bmatrix}\widetilde{P} & \mathsf{O}_{\nu} & \widetilde{Q} & \mathsf{O}_{\nu,q} & \mathsf{O}_{\nu,m}\end{bmatrix},\widetilde{P}\in\mathbb{S}^{\nu},\widetilde{Q}\in\mathbb{R}^{\nu\times d\nu}\\
&\boldsymbol{\Lambda}=\begin{bmatrix}P & Q\end{bmatrix},\ \widetilde{\boldsymbol{\Lambda}}=\begin{bmatrix}\widetilde{P} & \widetilde{Q}\end{bmatrix},\ \mathbf{N}=\mathbf{B}\boldsymbol{\mathfrak{K}},\ \ \widetilde{\mathbf{N}}=\mathbf{B}\widetilde{\boldsymbol{\mathfrak{K}}} \qquad (19)\\
&\widetilde{\boldsymbol{\mathfrak{K}}}=\begin{bmatrix}\widetilde{K}_1 & \widetilde{K}_2 & \widetilde{K}_3 & \mathsf{O}_{p,(q+m)}\end{bmatrix}\in\mathbb{R}^{p\times(2\nu+d\nu+q+m)}
\end{aligned}$$

in (16) with $\ell=2\nu+d\nu+q+m$ and $Z\oplus(I_\nu-Z)\succ 0$.

Applying inequality (17) with (19) to $\boldsymbol{\Theta}$ in (18) produces

$$\begin{aligned}
\boldsymbol{\Theta}&=\widehat{\boldsymbol{\Phi}}+\mathsf{Sy}\left(\mathbf{P}^\top\mathbf{B}\boldsymbol{\mathfrak{K}}\right)\preceq\widehat{\boldsymbol{\Phi}}+\mathsf{Sy}\left(\widetilde{\mathbf{P}}^\top\mathbf{N}+\mathbf{P}^\top\widetilde{\mathbf{N}}-\widetilde{\mathbf{P}}^\top\widetilde{\mathbf{N}}\right)\\
&\quad+\begin{bmatrix}\mathbf{P}^\top-\widetilde{\mathbf{P}}^\top & \mathbf{N}^\top-\widetilde{\mathbf{N}}^\top\end{bmatrix}\left[Z\oplus(I_\nu-Z)\right]^{-1}[*]\prec 0, \quad (20)
\end{aligned}$$

where we enforced that the upper bound of $\boldsymbol{\Theta}$ is negative-definite. Furthermore, $\mathbb{S}^{\nu}\ni Z\succ 0$ is a new decision variable. Applying the Schur complement [40] to (20) yields that (20) holds if and only if

$$\begin{bmatrix}\widehat{\boldsymbol{\Phi}}+\mathsf{Sy}\left(\widetilde{\mathbf{P}}^\top\mathbf{N}+\mathbf{P}^\top\widetilde{\mathbf{N}}-\widetilde{\mathbf{P}}^\top\widetilde{\mathbf{N}}\right) & \mathbf{P}^\top-\widetilde{\mathbf{P}}^\top & \mathbf{N}^\top-\widetilde{\mathbf{N}}^\top\\ * & -Z & \mathsf{O}_\nu\\ * & * & Z-I_\nu\end{bmatrix}\prec 0. \quad (21)$$

It appears that inequality (21) can be handled by standard SDP solvers [33] provided that the values of $\widetilde{\mathbf{P}}$ and $\widetilde{\boldsymbol{\mathfrak{K}}}$ are given.

By integrating all the preceding steps, we formulate Algorithm 1 in accordance with the expositions in [39], where vector $\mathbf{x}$ contains all the unknowns in $Z$ and Theorem 1, and $\|\mathbf{x}\|_\infty:=\max_i x_i$ and $\mathbf{v}^\top(\boldsymbol{\Lambda},\boldsymbol{\mathfrak{K}}):=\begin{bmatrix}\mathbf{vec}^\top(\boldsymbol{\Lambda}) & \mathbf{vec}^\top(\boldsymbol{\mathfrak{K}}))\end{bmatrix}$. Furthermore, $\rho_1,\rho_2>0$ and $\varepsilon>0$ are given constants for regularization and setting the error tolerance, respectively. To initialize Algorithm 1, we can simply substitute the values in (8) that can always be obtained via an explicit construction of PC (1) for IDS $\dot{\boldsymbol{x}}(t)=A\boldsymbol{x}(t)+B\boldsymbol{u}(t-r)$.

---

**Algorithm 1:** Inner Convex Approximation for Theorem 1

---

**begin**

  **solve** Theorem 1 with (8) **return** $\boldsymbol{\Lambda}=\begin{bmatrix}P & Q\end{bmatrix}$

  **solve** Theorem 1 with preceding $\boldsymbol{\Lambda}=\begin{bmatrix}P & Q\end{bmatrix}$ **return** $\boldsymbol{\mathfrak{K}}$.

  **update** $\widetilde{\boldsymbol{\Lambda}}\longleftarrow\boldsymbol{\Lambda}$, $\widetilde{\boldsymbol{\mathfrak{K}}}\longleftarrow\boldsymbol{\mathfrak{K}}$,

  **solve** $\min_{\mathbf{x}}\mathsf{tr}\left[\rho_1[*]\left(\boldsymbol{\Lambda}-\widetilde{\boldsymbol{\Lambda}}\right)+\rho_2[*]\left(\boldsymbol{\mathfrak{K}}-\widetilde{\boldsymbol{\mathfrak{K}}}\right)\right]$ subject to constraints (14b), (19) and (21), **return** $\boldsymbol{\Lambda}$ and $\boldsymbol{\mathfrak{K}}$

  **while** $\dfrac{\left\|\mathbf{v}(\boldsymbol{\Lambda},\boldsymbol{\mathfrak{K}})-\mathbf{v}(\widetilde{\boldsymbol{\Lambda}},\widetilde{\boldsymbol{\mathfrak{K}}})\right\|_\infty}{\left\|\mathbf{v}(\widetilde{\boldsymbol{\Lambda}},\widetilde{\boldsymbol{\mathfrak{K}}})\right\|_\infty+1}\geq\varepsilon$ **do**

    **update** $\widetilde{\boldsymbol{\Lambda}}\longleftarrow\boldsymbol{\Lambda}$, $\widetilde{\boldsymbol{\mathfrak{K}}}\longleftarrow\boldsymbol{\mathfrak{K}}$;

    **solve** $\min_{\mathbf{x}}\mathsf{tr}\left[\rho_1[*]\left(\boldsymbol{\Lambda}-\widetilde{\boldsymbol{\Lambda}}\right)+\rho_2[*]\left(\boldsymbol{\mathfrak{K}}-\widetilde{\boldsymbol{\mathfrak{K}}}\right)\right]$ subject to (14b), (19) and (21), **return** $\boldsymbol{\Lambda}$ and $\boldsymbol{\mathfrak{K}}$;

  **end**

**end**

---

## IV. Numerical Examples and Simulations

Computations were carried out via the SDP solvers Mosek [33] and SDPT3 [34], and all our SDP programs were coded using the optimization parser Yalmip [47] in ©MATLAB 2023b.

Consider (4) and (13) with $r=5$ and state space parameters

$$\begin{aligned}
&A=\begin{bmatrix}-1 & 1\\ 0 & 0.1\end{bmatrix},\ B=\begin{bmatrix}\beta(t)\\ 1+\beta(t)\end{bmatrix},\ D_1=\begin{bmatrix}0.1\\ -0.1\end{bmatrix},\\
&C_1=\begin{bmatrix}-0.3 & 0.4 & 0.1\\ -0.3 & 0.1 & -0.1\end{bmatrix},\ C_2=\begin{bmatrix}0 & 0.2 & 0\\ -0.2 & 0.1 & 0\end{bmatrix},\\
&\widetilde{C}_3(\tau)=\begin{bmatrix}0.2+0.1\mathsf{e}^{\tau} & 0.1 & 0.12\mathsf{e}^{3\tau}\\ -0.2 & 0.3+0.14\mathsf{e}^{2\tau} & 0.11\mathsf{e}^{3\tau}\end{bmatrix}, \qquad (22)\\
&D_2=0.12,\ \ D_3=\begin{bmatrix}0.14\\ 0.1\end{bmatrix},\\
&J_1=-\gamma I_2,\ \ \widetilde{J}=I_2,\ \ J_2=\mathbf{0}_2,\ \ J_3=\gamma>0
\end{aligned}$$

with $n=2$, $q=p=1$, $\nu=3$, where $A$ is non-Hurwitz, and $\beta(\cdot)$ satisfies $\widetilde{\forall}t\in\mathbb{R}$, $\beta(t)=0$, and $\gamma$ represents the $\mathscr{L}^2$ gain objective. Function $\beta(t)$ can represent glitches or noise that have countable nonzero values affecting input gain matrix $B$, which could not be incorporated if we were to use the traditional derivative for $\dot{\boldsymbol{x}}(t)$ in (4). As our methods are formulated for FDEs in the extended sense, $\beta(t)$ is regarded as 0 in our analysis w.r.t the Lebesgue measure. In fact, $\beta(t)$ can be added to any matrix in (22) without changing the calculus. This serves as a representative example showing the advantages of the proposed framework in modeling delay systems.

*To the best of our knowledge, no existing disturbance rejection methods can address (4) with the above parameters due to the presence of $\widetilde{C}_3(\tau),\beta(t)$ with $B\neq D_1$, where the columns of $B,D_1$ are linearly independent. Among the disturbance rejection methods in [12]–[15], quantitative tuning mechanisms have not been provided for the controller gains values which may improve the performance of rejection, as the methods require known gains.*

We first need to prescribe an initial value for the controller gain $\boldsymbol{\mathfrak{K}}$ on the basis of (3) and (1). Let $X=-1$ and $K=\begin{bmatrix}0 & -2\end{bmatrix}$ for the parameters in (8), which ensures that $A+BK$ is Hurwitz. By the structures of $(KA-XK)\,\mathsf{e}^{-A\tau}B=-2.2\mathsf{e}^{-0.1\tau}$ and the functions in $\widetilde{C}_3(\tau)$, we can employ

$$\boldsymbol{f}(\tau)=\begin{bmatrix}\mathop{\mathsf{Row}}_{i=0}^{\delta}\tau^i & \mathsf{e}^{\tau} & \mathsf{e}^{2\tau} & \mathsf{e}^{3\tau} & \mathsf{e}^{-0.1\tau}\end{bmatrix}^\top,\quad \delta\in\mathbb{N}\cup\{0\} \quad (23)$$

that is associated with $M=\begin{bmatrix}\mathbf{0}_\delta^\top & 0\\ \mathop{\mathsf{diag}}_{j=1}^{\delta} j & \mathbf{0}_\delta\end{bmatrix}\oplus 1\oplus 2\oplus 3\oplus(-0.1)$ for Assumption 1 to equivalently represent all the DDs in (22). With (23) and controller (7), we can immediately construct

$$\begin{aligned}
C_3=&\left[\begin{matrix}0.2 & 0.1 & 0 & \mathbf{0}_{\delta\nu}^\top & 0.1 & 0 & 0 & 0 & 0 & 0\\ -0.2 & 0.3 & 0 & \mathbf{0}_{\delta\nu}^\top & 0 & 0 & 0 & 0 & 0.14 & 0\end{matrix}\right.\\
&\qquad\left.\begin{matrix}0 & 0 & 0.12 & 0 & 0 & 0\\ 0 & 0 & 0.11 & 0 & 0 & 0\end{matrix}\right]\left(\sqrt{\mathfrak{F}}\otimes I_\nu\right), \qquad (24)\\
K_1=&\begin{bmatrix}0 & -3.6272 & -3\end{bmatrix},\quad K_2=\mathbf{0}_\nu^\top,\\
K_3=&\ \Gamma=\begin{bmatrix}\mathbf{0}_{\nu+\delta\nu}^\top & \mathbf{0}_{3\nu}^\top & 0 & 0 & -2.2\end{bmatrix}\left(\sqrt{\mathfrak{F}}\otimes I_\nu\right)
\end{aligned}$$

for (8) and Assumption 1, where $\sqrt{\mathfrak{F}}$ and $\sqrt{\mathfrak{F}^{-1}}$ are computed by the `vpaintegral` and `sqrtm` functions with variable precision in ©MATLAB to ensure accuracy ($\sqrt{\mathfrak{F}^{-1}}$ computing time $\approx 0.23$s).

We next applied Theorem 1 to (22) with (24) and $\delta=0$, yielding feasible solutions and showing that the resulting controller in the form of (3) ensures $\min\gamma=1.3915$. Note that $\mathbf{0}_0=[\,]_{0\times 1}$ and $\mathsf{diag}_{j=1}^{0}\,j=[\,]_{0\times 0}$ are empty matrices. Then, we computed a new $\boldsymbol{\mathfrak{K}}$ that guarantees $\gamma=1.3452$ by plugging in the resulting $P$, $Q$ into Theorem 1 again. With $\mathfrak{K}$ and $P, Q$ available, Algorithm 1 can be initialized. The results produced by Algorithm 1, with $\rho_1 = \rho_2 = 0.001$, are summarized in Table I, where the spectra abscissae of the resulting nominal CLSs were computed by the method in [48] via the Matlab toolbox `https://cdlab.uniud.it/software#eigAM-eigTMN`.

TABLE I: $\min \gamma$ COMPUTED BY ALGORITHM 1 WITH $\delta = 0$

| $\min \gamma$ | 1.1796 | 1.0984 | 1.0259 | 0.9614 |
|---|---|---|---|---|
| Number of Iterations | 100 | 200 | 300 | 400 |
| Spectra Abscissa | −0.1783 | −0.1501 | −0.1888 | −0.2189 |

Due to limited space, we only present the resulting $K_i$ after 400 iterations in (25). A crucial note for the readers is that we defined $F(\tau) = \sqrt{\mathfrak{F}^{-1}} \boldsymbol{f}(\tau) \otimes I_\nu$ in Assumption 1 and (7) and (29), where $\boldsymbol{f}(\tau)$ is normalized via $\sqrt{\mathfrak{F}^{-1}}$. An **incorrect** expression $K_3(\boldsymbol{f}(\tau) \otimes I_\nu)$ should never be used to obtain $K_3 F(\tau)$ in (7) for the resulting controller, as this will certainly cause errors.

$$
\begin{aligned}
K_1 &= \begin{bmatrix} 0.0776 & -0.3801 & -0.3994 \end{bmatrix}, \\
K_2 &= \begin{bmatrix} 0.0275 & -0.0317 & 0.012 \end{bmatrix}, \\
K_3 &= \Big[ -0.0250 \quad 0.0435 \quad -0.4908 \quad 0.0173 \quad \cdots \\
&\qquad 0.0312 \quad -0.2594 \quad 0.0262 \quad 0.0171 \quad \cdots \\
&\qquad\quad -0.1429 \quad 0.0085 \quad 0.0146 \quad -0.1085 \quad \cdots \\
&\qquad\qquad -0.0443 \quad 0.0142 \quad -0.4931 \Big].
\end{aligned} \tag{25}
$$

We repeated the above procedures but with $\delta = 3$ for $\boldsymbol{f}(\cdot)$ in (23). Every $\min \gamma$ values (See Table II) guaranteed by the resulting controllers are smaller than those in Table I under the same numbers of iterations. This shows that adding more functions to $\boldsymbol{f}(\cdot)$ in our controller (7) can improve the feasibility of Theorem 1 and Algorithm 1, which results in controllers with better performance.

TABLE II: $\min \gamma$ COMPUTED BY ALGORITHM 1 WITH $\delta = 3$

| $\min \gamma$ | 1.1783 | 1.0968 | 1.0226 | 0.9555 |
|---|---|---|---|---|
| Number of Iterations | 100 | 200 | 300 | 400 |
| Spectral Abscissa | −0.1734 | −0.1892 | −0.2177 | −0.1163 |

Clearly, Tables I-II show that more iterations can lead to smaller $\min \gamma$ values but with increased computational complexities. The results also show the advantage of the LDSFs in (7), which can ensure better performance than the controller structure in (3).

The SPA values in Tables I-II are all negative, confirming that all the resulting CLSs are exponentially stable with $\boldsymbol{w}(t) \equiv \mathbf{0}_q$. Moreover, controller (7) satisfies the properties in Conclusion 1 as all CLSs are of the retarded type. This ensures all the DDs of the CLSs can be safely implemented using ordinary quadratures [31] without causing instability problems, provided that the numerical accuracy reaches certain level [31, Th. 2].

For numerical simulations, we select the resulting controller parameters in (25) for (7), corresponding to $\min \gamma = 0.9614$ in Table I. Let $t_0 = 0$, $\boldsymbol{z}(t) = \mathbf{0}_2, \forall t < 0$, and $\boldsymbol{\phi}(\tau) = \begin{bmatrix} 1 & 2 & 0.5 \end{bmatrix}^\top, \forall \tau \in [-5, 0]$ be the initial condition, and $w(t) = 10 \sin(10t) t^2 \mathrm{e}^{-0.8t}$ be a scalar disturbance. We employed the Band-Limited White Noise block in Simulink, with `Sample time` = 0.002 and default values for `Seed` and `Noise power`, to generate a white noise signal $n(t)$ for $\beta(t) = n(t) \left( \mathit{1}(t) - \mathit{1}(t - 10) \right)$ in (22), where $\mathit{1}(t)$ is the Heaviside step function. Since $\forall t \geq 10, \beta(t) = 0$ holds true and $\beta(t)$ is sampled into a discrete sequence, function $\beta(t)$ has only a **finite** number of nonzero data points during simulation, which satisfies $\widetilde{\forall} t \geq 0$, $\beta(t) = 0$ mathematically, as in (22). Moreover, all the DDs in the CLS were discretized via the standard trapezoidal rule with 500 sample points. The simulation was performed in Simulink using the `ode8` solver with $0.002s$ as the fundamental sampling time. The results in Figures 2–3 clearly show that the goal of $\mathscr{L}^2$ stabilization is achieved by the resulting controller for the CLS, while the effects of noise were gradually suppressed and eventually extinguished.

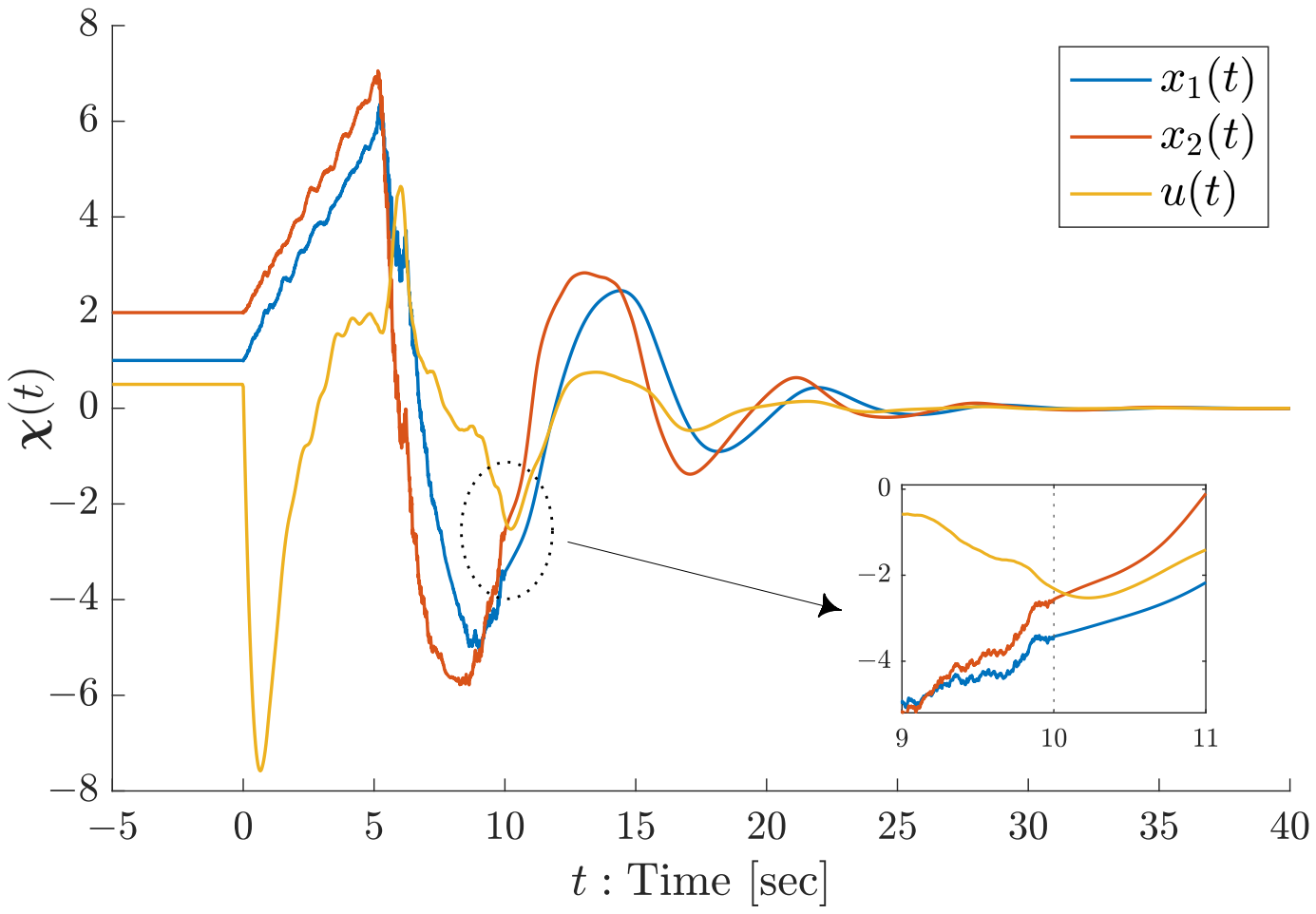

Fig. 2: Plots of each CLS state $\boldsymbol{\chi}(t) = \begin{bmatrix} \boldsymbol{x}^\top(t) & u^\top(t) \end{bmatrix}^\top$.

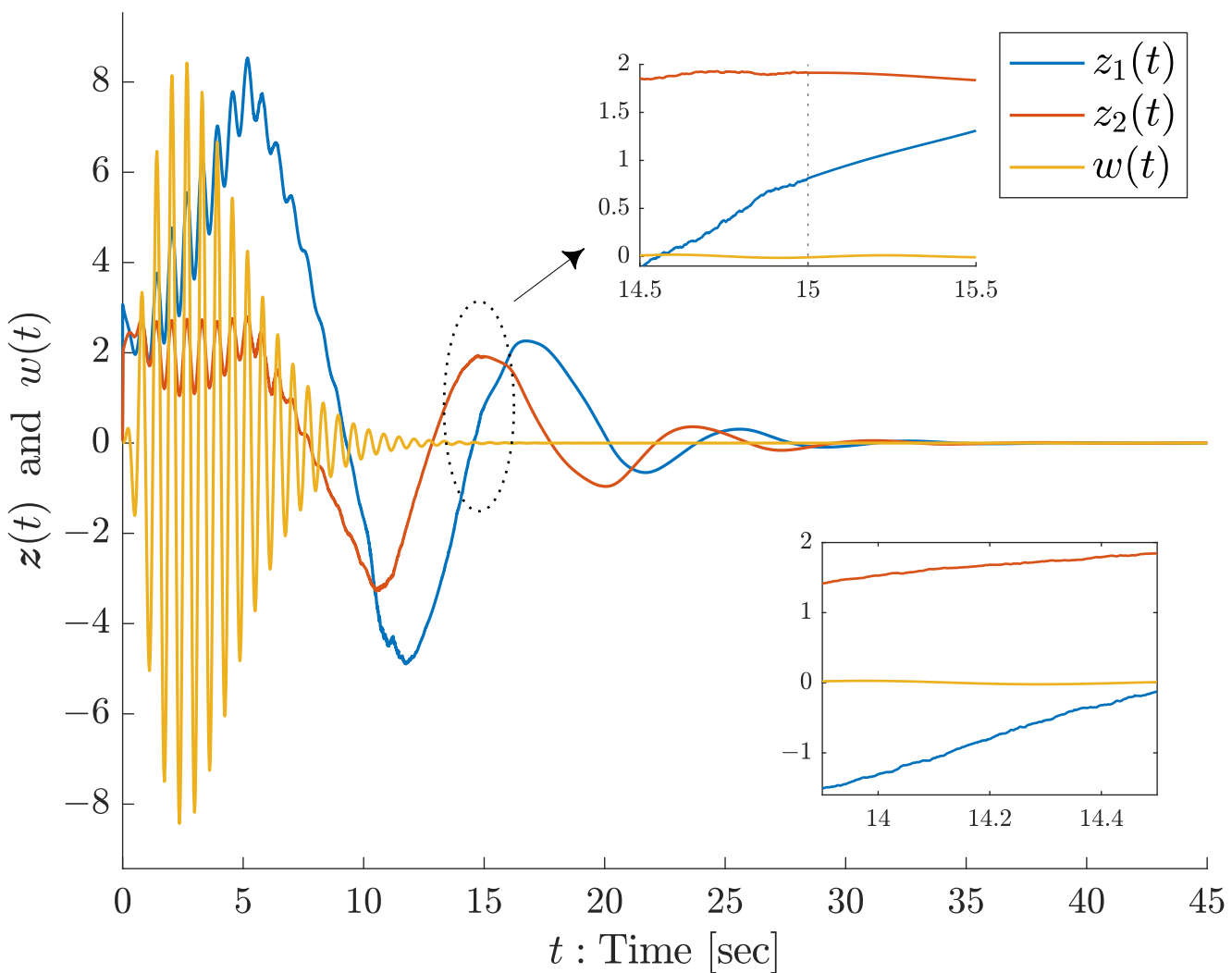

Fig. 3: Plots of regulated output $\boldsymbol{z}(t)$ and disturance $w(t)$.

## V. CONCLUSION AND FUTURES STUDIES

Our proposed framework offers an original perspective on addressing the input delay stabilization problems with dissipative constraints characterized by the SRF in (13). The key components of our methodology include the use of the FDEs in the extended sense, the application of the KD in Assumption 1, a novel concept of LDSFs in (7), the construction of KF (29) connected to the integral kernels in Assumption 1 and (7), and the inner convex approximation in Algorithm 1. Moreover, our design also takes into account potential disturbances caused by the operating environment as shown in Figure 1, thus the resulting controller (7) itself is also robust. The initial values of Algorithm 1 can always be supplied by the gains of (3), and the algorithm's feasibility can be enhanced by adding more functions to $\boldsymbol{f}(\cdot)$ in (29) or (7) satisfying the conditions in (6). Since no products exist in (10) between the controller gains in (7) and the state space parameters in (4), we will further investigate the advantages of our LDSF concept in Part II of this series in dealing with uncertainties. We will illustrate how the structure of our LDSF can facilitate the design of resilient (non-fragile) controllers for uncertain IDSs with dissipative constraints, by utilizing the lemma we have developed for the analysis of linear fractional uncertainties.

## Appendix

### A. Integral Inequality [49, Th. 1]

We define the weighted Lebesgue function space

$$\mathscr{L}^2_{\varpi}\left(\mathcal{K};\mathbb{R}^d\right) := \left\{\boldsymbol{\varphi}(\cdot) \in \mathcal{M}\left(\mathcal{K};\mathbb{R}^d\right) : \left\|\boldsymbol{\varphi}(\cdot)\right\|_{\varpi} < \infty\right\} \quad (26)$$

with seminorm $\|\boldsymbol{\varphi}(\cdot)\|^2_{\varpi} := \int_{\mathcal{K}} \varpi(\tau)\boldsymbol{\varphi}^\top(\tau)\boldsymbol{\varphi}(\tau)\mathsf{d}\tau$, where the weight function $\varpi(\cdot) \in \mathcal{M}(\mathcal{K};\mathbb{R}_{\geq 0})$ has only countably infinite or finite zero values. Furthermore, $\mathcal{K}$ is Lebesgue measurable [38] and satisfies $\mathcal{K} \subseteq \mathbb{R} \cup \{\pm\infty\}$ and $\int_{\mathcal{K}} \mathsf{d}\tau \neq 0$.

*Theorem 2:* Let $\varpi(\cdot)$, $\mathcal{K}$ and $d \in \mathbb{N}$ in (26) be given, and suppose that $\boldsymbol{g}(\cdot) \in \mathscr{L}^2_{\varpi}\left(\mathcal{K};\mathbb{R}^d\right)$ satisfies [40, Th. 7.2.10] inequality

$$\int_{\mathcal{K}} \varpi(\tau)\boldsymbol{g}(\tau)\boldsymbol{g}^\top(\tau)\mathsf{d}\tau \succ 0. \quad (27)$$

Then, inequality

$$\int_{\mathcal{K}} \varpi(\tau)\boldsymbol{y}^\top(\tau)Y\boldsymbol{y}(\tau)\mathsf{d}\tau \geq [*]\left(\mathfrak{G}^{-1}\otimes Y\right)\int_{\mathcal{K}} \varpi(\tau)G(\tau)\boldsymbol{y}(\tau)\mathsf{d}\tau \quad (28)$$

holds for all $\boldsymbol{y}(\cdot) \in \mathscr{L}^2_{\varpi}(\mathcal{K};\mathbb{R}^\nu)$ and $Y \succ 0$, where $G(\tau) = \boldsymbol{g}(\tau)\otimes I_\nu$ and $\mathfrak{G}^{-1} = \int_{\mathcal{K}} \varpi(\tau)\boldsymbol{g}(\tau)\boldsymbol{g}^\top(\tau)\mathsf{d}\tau$. Moreover, the optimality of (28) is guaranteed by the least-square principle [49, Cor. 1].

The above inequality generalizes the Bessel Legendre inequality in [44], [45] and its derivatives.

### B. Proof of Theorem 1

The proof is established by the construction of KF

$$\begin{aligned}\mathsf{v}(\boldsymbol{\chi}_t(\cdot)) := [*]\begin{bmatrix} P & Q \\ * & R\end{bmatrix}\begin{bmatrix}\boldsymbol{\chi}(t) \\ \int_{-r}^0 F(\tau)\boldsymbol{\chi}(t+\tau)\mathsf{d}\tau\end{bmatrix} \\ + \int_{-r}^0 \boldsymbol{\chi}^\top(t+\tau)\left[S+(\tau+r)U\right]\boldsymbol{\chi}(t+\tau)\mathsf{d}\tau\end{aligned} \quad (29)$$

with $F(\tau), \mathfrak{F}$ in Assumption 1, where $F(\tau) = \sqrt{\mathfrak{F}^{-1}}\boldsymbol{f}(\tau)\otimes I_\nu$. First, we note that $\mathsf{s}(\boldsymbol{z}(t),\boldsymbol{w}(t))$ in (13) can be rewritten as

$$\begin{aligned}\mathsf{s}(\boldsymbol{z}(t),\boldsymbol{w}(t)) = \boldsymbol{z}^\top(t)\widetilde{J}^\top J_1^{-1}\widetilde{J}\boldsymbol{z}(t) + \mathsf{Sy}\left[\boldsymbol{z}^\top(t)J_2\boldsymbol{w}(t)\right] \\ + \boldsymbol{w}^\top(t)J_3\boldsymbol{w}(t).\end{aligned} \quad (30a)$$

Following from the expression of $\boldsymbol{z}(t)$ in (10), we have

$$\boldsymbol{z}^\top(t)\widetilde{J}^\top J_1^{-1}\widetilde{J}\boldsymbol{z}(t) = \boldsymbol{\vartheta}^\top(t)\boldsymbol{\Sigma}^\top\widetilde{J}^\top J_1^{-1}\widetilde{J}\boldsymbol{\Sigma}\boldsymbol{\vartheta}(t) \quad (30b)$$

via the matrix $\boldsymbol{\Sigma}$ and vectored-value function $\boldsymbol{\vartheta}(t)$ in (10).

Next, we differentiate (weak derivative) $\mathsf{v}(\boldsymbol{\chi}_t(\cdot))$ in (29) along the trajectory of CLS (10) and consider (30a)–(30b), and the identity

$$\begin{aligned}\frac{\mathsf{d}}{\mathsf{d}t}\int_{-r}^0 F(\tau)\boldsymbol{\chi}(t+\tau)\mathsf{d}\tau = F(0)\boldsymbol{\chi}(t) - F(-r)\boldsymbol{\chi}(t-r) \\ -(M\otimes I_\nu)\int_{-r}^0 F(\tau)\boldsymbol{\chi}(t+\tau)\mathsf{d}\tau = \mathbf{F}\boldsymbol{\vartheta}(t)\end{aligned} \quad (31)$$

derived from an application of (6) with $\mathbf{F}$ in (15d). Then we deduce

$$\begin{aligned}\widetilde{\forall} t \geq t_0, \quad \dot{\mathsf{v}}(\boldsymbol{\chi}_t(\cdot)) - \mathsf{s}(\boldsymbol{z}(t),\boldsymbol{w}(t)) = \\ \boldsymbol{\vartheta}^\top(t)\,\mathsf{Sy}\left(\begin{bmatrix} I_\nu & \mathsf{O}_{\nu,d\nu} \\ \mathsf{O}_\nu & \mathsf{O}_{\nu,d\nu} \\ \mathsf{O}_{d\nu,\nu} & I_{d\nu} \\ \mathsf{O}_{q,\nu} & \mathsf{O}_{q,d\nu}\end{bmatrix}\begin{bmatrix} P & Q \\ * & R\end{bmatrix}\begin{bmatrix}\mathbf{A}+\mathbf{BK} \\ \mathbf{F}\end{bmatrix}\right)\boldsymbol{\vartheta}(t) \\ + \boldsymbol{\chi}^\top(t)\left(S+rU\right)\boldsymbol{\chi}(t) - \boldsymbol{\chi}^\top(t-r)S\boldsymbol{\chi}(t-r) \\ - \int_{-r}^0 \boldsymbol{\chi}^\top(t+\tau)U\boldsymbol{\chi}(t+\tau)\mathsf{d}\tau - \boldsymbol{w}^\top(t)J_3\boldsymbol{w}(t) \\ - \boldsymbol{\vartheta}^\top(t)\,\mathsf{Sy}\left(\boldsymbol{\Sigma}^\top\widetilde{J}^\top J_1^{-1}\widetilde{J}\boldsymbol{\Sigma} + \begin{bmatrix}\mathsf{O}_{(2\nu+d\nu),m} \\ J_2^\top\end{bmatrix}\boldsymbol{\Sigma}\right)\boldsymbol{\vartheta}(t).\end{aligned} \quad (32)$$

Assume that $U \succ 0$, if follows from utilizing (28) on $\int_{-r}^0 [*]U\boldsymbol{\chi}(t+\tau)\mathsf{d}\tau$ in (32) with $\varpi(\tau) = 1$, $\boldsymbol{g}(\tau) = \sqrt{\mathfrak{F}^{-1}}\boldsymbol{f}(\tau)$ and $Y = U$, that

$$\begin{aligned}\int_{-r}^0 \boldsymbol{\chi}^\top(t+\tau)U\boldsymbol{\chi}(t+\tau)\mathsf{d}\tau \geq \\ \int_{-r}^0 \boldsymbol{\chi}^\top(t+\tau)F^\top(\tau)\mathsf{d}\tau\,(I_d\otimes U)\int_{-r}^0 F(\tau)\boldsymbol{\chi}(t+\tau)\mathsf{d}\tau\end{aligned} \quad (33)$$

with $F(\tau)$ in (29). Next, applying (33) to (32) yields inequality

$$\begin{aligned}\widetilde{\forall} t \geq t_0, \quad \dot{\mathsf{v}}(\boldsymbol{\chi}_t(\cdot)) - \mathsf{s}(\boldsymbol{z}(t),\boldsymbol{w}(t)) \\ \leq \boldsymbol{\vartheta}^\top(t)\left(\boldsymbol{\Psi} - \boldsymbol{\Sigma}^\top\widetilde{J}^\top J_1^{-1}\widetilde{J}\boldsymbol{\Sigma}\right)\boldsymbol{\vartheta}(t)\end{aligned} \quad (34a)$$

following from the structural components of $\boldsymbol{\vartheta}(t)$ in (10), where

$$\begin{aligned}\boldsymbol{\Psi} = \mathsf{Sy}\left(\begin{bmatrix} P & Q \\ \mathsf{O}_\nu & \mathsf{O}_{\nu,d\nu} \\ Q^\top & R \\ \mathsf{O}_{q,\nu} & \mathsf{O}_{q,d\nu}\end{bmatrix}\begin{bmatrix}\mathbf{A}+\mathbf{BK} \\ \mathbf{F}\end{bmatrix} + \begin{bmatrix}\mathsf{O}_{(2\nu+d\nu),m} \\ J_2^\top\end{bmatrix}\boldsymbol{\Sigma}\right) \\ - \left(\left[-S-rU\right]\oplus S\oplus(I_d\otimes U)\oplus J_3\right)\end{aligned} \quad (34b)$$

comprises all the matrices induced by $\dot{\mathsf{v}}(\boldsymbol{\chi}_t(\cdot)) - \mathsf{s}(\boldsymbol{z}(t),\boldsymbol{w}(t))$ in (32) excluding the terms in (30b). By examining the upper bound in (34a) with $\boldsymbol{\vartheta}(t)$ in (10), we find that (12) is true if $U \succ 0$ and

$$\boldsymbol{\Psi} - \boldsymbol{\Sigma}^\top\widetilde{J}^\top J_1^{-1}\widetilde{J}\boldsymbol{\Sigma} \prec 0. \quad (35)$$

Now applying the Schur complement to (35) shows that (14b)

$$\begin{bmatrix}\boldsymbol{\Psi} & \boldsymbol{\Sigma}^\top\widetilde{J}^\top \\ * & J_1\end{bmatrix} = \mathbf{P}^\top\begin{bmatrix}\boldsymbol{\Pi} & \mathsf{O}_{\nu,m}\end{bmatrix} + \boldsymbol{\Phi} \prec 0$$

holds if and only if (35) holds with $U \succ 0$, given $J_1^{-1} \prec 0$ in (13). Thus, KF (29) satisfies (12) if (14b) and $U \succ 0$ are feasible. In light of the structure of $\boldsymbol{\vartheta}(t)$ in (10) and $\boldsymbol{\Psi} \prec 0$ in (14b), we find that (11b) holds for some $\epsilon_3 > 0$ if (14b) and $U \succ 0$ are true.

Now we begin to prove that (29) satisfies (11a) if (14a) holds. First, note that $\forall\theta \in [-r,0]$, $\boldsymbol{\chi}_{t_0}(\theta) = \boldsymbol{\chi}(t_0+\theta) = \boldsymbol{\phi}(\theta)$ in (10). Let $t = t_0$ in (29) and $S \succ 0$; then, we find that $\exists\lambda > 0:$

$$\begin{aligned}\mathsf{v}(\boldsymbol{\phi}(\cdot)) \leq [*]\,\lambda\begin{bmatrix}\boldsymbol{\phi}(0) \\ \int_{-r}^0 F(\tau)\boldsymbol{\phi}(\tau)\mathsf{d}\tau\end{bmatrix} + \lambda\int_{-r}^0 \left\|\boldsymbol{\phi}(\tau)\right\|^2\mathsf{d}\tau \\ \leq [*]\lambda\boldsymbol{\phi}(0) + [*](I_d\otimes\lambda I_\nu)\int_{-r}^0 F(\tau)\boldsymbol{\phi}(\tau)\mathsf{d}\tau + r\lambda\left\|\boldsymbol{\phi}(\cdot)\right\|^2_\infty \\ \leq \boldsymbol{\phi}^\top(0)\lambda\boldsymbol{\phi}(0) + \int_{-r}^0 \boldsymbol{\phi}^\top(\tau)\lambda I_\nu\boldsymbol{\phi}(\tau)\mathsf{d}\tau + r\lambda\left\|\boldsymbol{\phi}(\cdot)\right\|^2_\infty \\ \leq (\lambda+2r\lambda)\left\|\boldsymbol{\phi}(\cdot)\right\|^2_\infty\end{aligned} \quad (36)$$

for any $\boldsymbol{\phi}(\cdot) \in \mathcal{C}\left([-r,0];\mathbb{R}^\nu\right)$, following from utilizing (28) with $\boldsymbol{g}(\tau) = \sqrt{\mathfrak{F}^{-1}}\boldsymbol{f}(\tau)$, $\varpi(\tau) = 1$, $Y = \lambda I_\nu$ and the property of quadratic form $\forall W \in \mathbb{S}^\nu, \exists\lambda > 0, \forall\mathbf{x}\setminus\{\mathbf{0}_\nu\}, \mathbf{x}^\top(\lambda I_\nu - W)\mathbf{x} > 0$. Thus, we have shown that the KF in (29) satisfies (11a) if (14a) holds.

By applying (28) again with $\boldsymbol{g}(\tau) = \sqrt{\mathfrak{F}^{-1}}\boldsymbol{f}(\tau)$, $\varpi(\tau) = 1$ and $Y = S$ to $\int_{-r}^0 \boldsymbol{\chi}^\top(t+\tau)S\boldsymbol{\chi}(t+\tau)\mathsf{d}\tau$ in (29) at $t = t_0$, we have

$$\int_{-r}^0 \boldsymbol{\phi}^\top(\tau)S\boldsymbol{\phi}(\tau)\mathsf{d}\tau \geq [*]\left(I_d\otimes S\right)\int_{-r}^0 F(\tau)\boldsymbol{\phi}(\tau)\mathsf{d}\tau \quad (37)$$

with $S \succ 0$. We then utilize (37) on (29) at $t = t_0$ finding that

$$\mathsf{v}(\boldsymbol{\phi}(\cdot)) \geq [*] \begin{bmatrix} P & Q \\ * & R + I_d \otimes S \end{bmatrix} \begin{bmatrix} \boldsymbol{\phi}(0) \\ \int_{-r}^{0} F(\tau)\boldsymbol{\phi}(\tau)\mathsf{d}\tau \end{bmatrix} + \int_{-r}^{0} (\tau + r)\boldsymbol{\phi}^{\top}(\tau) U \boldsymbol{\phi}(\tau)\mathsf{d}\tau \tag{38}$$

for any $\boldsymbol{\phi}(\cdot) \in \mathcal{C}([-r, 0]; \mathbb{R}^{\nu})$. Given (38) with the properties of positive definite matrices, we conclude that $\exists \epsilon_1 > 0 : \mathsf{v}(\boldsymbol{\phi}(\cdot)) \geq \epsilon_1 \|\boldsymbol{\phi}(0)\|^2$ for any $\boldsymbol{\phi}(\cdot)$ in (10) if (14a) holds.

Since we have shown that the KF in (29) satisfies (12) and (11b) if (14b) holds with $U \succ 0$, it has bee proved that the existence of feasible solutions to (14a) and (14b) indicates the existence of a KF in (29) satisfying dissipativity condition (12) and the stability criteria in (11a) and (11b).